\documentclass[11pt]{amsart} %

\usepackage{float}

\usepackage{epsfig}
\usepackage{epstopdf}
\usepackage{tikz}
\usepackage{array,amsmath, enumerate,  url, psfrag}
\usepackage{amssymb, amsaddr, fullpage}
\usepackage{graphicx,subfigure}
\usepackage{color}
\usepackage{amssymb}
\theoremstyle{plain}

\newtheorem{theorem}{Theorem}
\numberwithin{theorem}{section}

\numberwithin{lemma}{section}

\numberwithin{corollary}{section}

\newtheorem{conjecture}{Conjecture}
\numberwithin{conjecture}{section}
\newtheorem*{conjecture*}{Conjecture}

\newtheorem{observation}{Observation}
\numberwithin{observation}{section}

\title{Enhancing the Erd\H{o}s-Lov\'asz Tihany Conjecture for graphs with independence number two}

\author{Yue Wang$^{1}$,  Gexin Yu$^{2}$}

\address{
$^{1}$\small School of Mathematics, Shandong University, Jinan, Shandong, China.\\
$^{2}$\small Department of Mathematics, William \& Mary, Williamsburg, VA, USA.
}

\thanks{The work was done while the first author was at William \& Mary as a visiting student, partially  supported by the Chinese Scholarship Council.  The research of the last author was supported in part by a summer research grant from William \& Mary.}

\email{m15064013175@163.com(Y. Wang), gyu@wm.edu}

\date{\today}

\begin{document}

\maketitle

\begin{abstract}
Let $s\ge2$  and $t\ge2$ be integers. A graph $G$ is $(s,t)$-\emph{splittable} if $V(G)$ can be partitioned into two sets $S$ and $T$ such that $\chi(G[S])\geq s$ and $\chi(G[T])\geq t$. The well-known Erd\H{o}s-Lov\'asz Tihany Conjecture from 1968 states that every graph $G$ whose chromatic number $\chi(G)=s+t-1$ is more than its clique number $\omega(G)$ is $(s,t)$-splittable. In this paper, we prove an enhanced version of the Erd\H{o}s-Lov\'asz Tihany Conjecture for graphs with independence number two. That is, for every graph $G$ with $\chi(G)=s+t-1>\omega(G)+1$ is $(s,t+1)$-splittable. There are examples showing that this result is best possible.
\end{abstract}

\section{Introduction}\label{sect1}

All graphs considered in this paper are finite and without loops or multiple edges. Given a graph $G$, we write $n(G)$  for the number of vertices of
$G$, $\alpha(G)$ for its independence number, $\omega(G)$ for its clique number, $\chi(G)$ for its chromatic number, 
and $\overline{G}$ for the complement of $G$. Given a vertex set $A\subseteq V(G)$, the subgraph of $G$ induced by $A$, denoted $G[A]$, is the graph with vertex set $A$ and edge set $\{xy\in E(G):x,y\in A\}$.

Let $s\ge2$  and $t\ge2$ be integers. A graph $G$ is $(s,t)$-\emph{splittable} if $V(G)$ can be partitioned into two sets $S$ and $T$ such that $\chi(G[S])\geq s$ and $\chi(G[T])\geq t$.  In 1968, Erd\H{o}s and Lov\'asz~\cite{EL68} made the following famous conjecture:

\begin{conjecture}\label{con1}\emph{(Erd\H{o}s-Lov\'asz Tihany Conjecture)}.
For every graph $G$ with $\chi(G) > \omega(G)$ and any two integers $s,t\geq 2$ with $s+t =\chi(G)+1$, $G$ is $(s,t)$-splittable.
\end{conjecture}

The only settled cases of this conjecture are $(s,t)\in\{(2,2),(2,3),(2,4),(3,3),(3,4),(3,5)\},$  see \cite{BJ99, M86, S87a, S87b}. This conjecture is also known to be true for some special classes of graphs, such as line graphs of multigraphs (Kostochka and Stiebitz \cite{KS08}), quasi-line graphs and graphs with independence number two (Balogh, Kostochka, Prince and Stiebitz~\cite{BKPS09}).

A connected graph $G$ is {\em double-critical} if $\chi(G)=t$ but $\chi(G\backslash\{x,y\})=t-2$ for every edge $xy\in E(G)$.  The following well-known conjecture is the case of $s=2$ of Conjecture~\ref{con1}.

\begin{conjecture}\label{con2}\emph{(Double-Critical Graph Conjecture~\cite{EL68})}.
For $t\ge 3$, the only double-critical $t$-chromatic graph is $K_t$.
\end{conjecture}

From \cite{S87a}, Conjecture~\ref{con2} holds when $t\le 5$. For $t\ge 6$,  Conjecture \ref{con2} remains wide open, and we even do not know if every double-critical $t$-chromatic graph contains $K_4$ as a subgraph.

As Conjecture~\ref{con1} implies Conjecture~\ref{con2} which seems hopeless to prove at this moment, we would like to study a version of Conjecture~\ref{con1} that does not imply Conjecture~\ref{con2}.  In other words, is Conjecture~\ref{con1} true if $s,t\ge 3$?

The Conjecture~\ref{con1} can be greatly enhanced in line graphs \cite{WY}:  when $s,t\ge 3.5\ell+2$, for the line graph $G$ of some multigraph, $G$ is $(s,t+\ell)$-splittable. In this paper, we consider the graphs with independence number two and get the following result.

\begin{theorem} \label{th1}
Let $s$ and $t$ be arbitrary integers with $t\ge s\ge 2$. If a graph $G$ with $\alpha(G)=2$ and $\chi(G)=s+t-1> \omega(G)+1$, then  $G$ is $(s,t+1)$-splittable.
\end{theorem}

It is worth mentioning that there are examples showing that Theorem \ref{th1} is best possible.\\

\noindent\textbf{Example 1.} Let $s$ and $t$ be arbitrary integers with $t\ge s\ge 2$. There is a graph $G$ with $\alpha(G)=2$ and $\chi(G)=s+t-1> \omega(G)$ but $G$ is not $(s,t+1)$-splittable.

\begin{proof}
Let $X$ be a $K_{s+t-4}$ and $Y$ be a $C_5$. The graph $G$ can be partitioned into two parts $X$ and $Y$ so that every vertex in $X$ is adjacent to every vertex in $Y$. As $C_5$ has no triangle, it is obvious that $\omega(G)=s+t-4+2=s+t-2$. Note that $\overline{G}$ has no triangle, so $\alpha(G)=2$. Since every vertex in $X$ is adjacent to every vertex in $Y$, $\chi(G)=\chi(X)+\chi(Y)=s+t-1$. Clearly, $|V(G)|=|X|+|Y|=s+t+1$.

Suppose that $G$ is $(s,t+1)$-splittable. Then there is a partition $X_1,X_2$ of $G$ with $|X_1|=s$ and $|X_2|=t+1$ such that different vertices must belong to different color classes. Clearly, there must exist a part, say $X_1$ containing at least three vertices of $Y$. But these three vertices of $Y$ does not form a triangle and so two vertices belong to the same color class, a contradiction.
\end{proof}

\noindent\textbf{Example 2.} Let $s$ and $t$ be arbitrary integers with $t\ge s\ge 2$. There is a graph $G$ with $\alpha(G)=2$ and $\chi(G)=s+t-1> \omega(G)+1$ but $G$ is not $(s,t+2)$-splittable.

\begin{proof}
Let $X=K_{s+t-7}$, $Y=Z=C_5$. The graph $G$ can be partitioned into three parts $X, Y$ and $Z$ so that any two vertices in different parts are adjacent. As $C_5$ has no triangle, it is obvious that $\omega(G)=s+t-7+2+2=s+t-3$. And note that $\overline{G}$ has no triangle, so $\alpha(G)=2$. Since any two vertices in different parts are adjacent, $\chi(G)=\chi(X)+\chi(Y)+\chi(Z)=s+t-1$. Clearly, $|V(G)|=|X|+|Y|+|Z|=s+t+3$.

Suppose that $G$ is $(s,t+2)$-splittable. Then there is a partition $X_1,X_2$ of $G$ with $|X_1|=s$ and $|X_2|=t+2$ such that there is at most one color class containing two vertices. Clearly, there must exist a part containing at least three vertices of $C_5$. And these three vertices of $C_5$ does not form a triangle and so two vertices belong to the same color class. As $Y=Z=C_5$, there are two color classes containing two vertices, a contradiction.
\end{proof}

In the next section, we prove the main results.

\section{Proof of Theorem \ref{th1}}\label{sect2}

\begin{proof}
Let $G$ be a counterexample to the theorem and $|V(G)|=n$. We denote by $o(G)$ the number of odd components in the graph $G$.
\begin{observation}\label{ob1}\cite{BKPS09}
If $G$ is a graph with independence number 2, then $$\chi(G)=\max\limits_{P\subseteq V(G)}\left\{\frac{n(G)+o(\overline{G}-P)-|P|}2\right\}.$$
\end{observation}

According to Observation \ref{ob1}, there is a $P\subseteq V(G)$ such that
\begin{equation}
\chi(G)=\frac{n(G)+o(\overline{G}-P)-|P|}2.\label{eq2}
\end{equation}
We choose a largest such set $P$. We claim that every component of $\overline{G}-P$ is odd. Suppose to the contrary that there is a component $H'$ of $\overline{G}-P$ that is even. Choose a vertex $v\in V(H')$. Let $P'=P\cup \{v\}$. As $H'-v$ is odd, we get that $o(\overline{G}-P')\ge o(\overline{G}-P)+1$. It follows that $$\frac{n(G)+o(\overline{G}-P')-|P'|}2\ge\frac{n(G)+o(\overline{G}-P)+1-(|P|+1)}2=\chi(G),$$
a contradiction to the maximality of $P$.

We denote by $\mathcal{H}$ the set of all components of $\overline{G}-P$. Let $\mathcal{H}_i$ denote the set of all components of size $i$ in $\mathcal{H}$ and let $|\mathcal{H}_i|=k_i$. 

\noindent\textbf{Case 1.} $\omega(G)\ge s$.

Let $S_0$ be a set of $s$ vertices forming a clique in $G$ and $T_0=V(G)\backslash S_0$. Since $\alpha(G)=2$, we have that $$t\ge\chi(G[T_0])\ge\frac{n-|S_0|}{2}=\frac{n-s}2.$$
Adding $s$ to both sides, we have that
\begin{equation}
\chi(G)+1=s+t\ge\frac{n+s}2.\label{eq1}
\end{equation}
Combining \eqref{eq1} and \eqref{eq2}, we have that
\begin{equation}
o(\overline{G}-P)\ge s-2+|P|.\label{eq3}
\end{equation}

\textbf{Subcase 1.1} $o(\overline{G}-P)-k_1\ge2$.
Since $\alpha(G)=2$, $\overline{G}$ is triangle-free. It follows that for $i\ge3$, each component in $\mathcal{H}_i$ contains a pair $\{x,y\}$ of non-adjacent vertices. As $o(\overline{G}-P)-k_1\ge2$, we can choose at least two such pairs of non-adjacent vertices $\{x_1,y_1\}\subseteq H_i$, $\{x_2,y_2\}\subseteq H_j$ and $H_i\neq H_j$, where $H_i,H_j\in \mathcal{H}-\mathcal{H}_1$. By \eqref{eq3}, we can choose a set $S'$ of $s-4$ vertices, each from a different component of $\mathcal{H}-H_i-H_j$. Let $S=\{x_1,y_1,x_2,y_2\}\cup S'$. According to the construction of $S$, $S$ induces an $s$-clique in $G$. Since $|V(H_i)-\{x_1,y_1\}|$ and $|V(H_j)-\{x_2,y_2\}|$ are odd, $o(\overline{G}-S-P)\ge o(\overline{G}-P)-(s-4)$. By Observation \ref{ob1},
\begin{align*}
\chi(G-S)&\ge \frac{(n-s)+o(\overline{G}-S-P)-|P|}2\\
&\ge\frac{n-s+o(\overline{G}-P)-(s-4)-|P|}2\\
&=\chi(G)-s+2=t+1.
\end{align*}
This contradicts our assumption.

\textbf{Subcase 1.2} $o(\overline{G}-P)-k_1<2$. We claim that $o(\overline{G}-P)\neq k_1$. Suppose to the contrary that $o(\overline{G}-P)= k_1$. 
By definition, $o(\overline{G}-P)=n-|P|$. By \eqref{eq2}, $\omega(G)\ge k_1=\chi(G)$, a contradiction. We denote by $H_i\in \mathcal{H}_i$ where $i=1$ (mod 2) and $i\ge3$. 
So $o(\overline{G}-P)= k_1+1$ and $n-|P|= k_1+i$.  It follows that $\chi(G)= k_1+\frac{1+i}2$ from \eqref{eq2}. We also know that $\omega(G)\ge k_1+\alpha(H_i)$. As $\omega(G)\le\chi(G)-2$, we have that
\begin{equation}
\alpha(H_i)\le\frac{i-3}2.\label{eq1.2}
\end{equation}
As $\alpha(H_i)\ge2$, $i\ge7$ by \eqref{eq1.2}. Since $i>R(3,3)=6$ and $\overline{G}$ has no triangle, it implies $\alpha(H_i)\ge3$. Thus, we finally get that $i\ge9$ by \eqref{eq1.2}. It follows that $\alpha(H_i)\ge4$ from the fact that $R(3,4)=9$, for $i\ge9$.
Choose an independent set $\{x_1,x_2,x_3,x_4\}\subseteq H_i$. By \eqref{eq3}, we can choose a set $S'$ of $s-4$ vertices, each from a different component of $\mathcal{H}-H_i$. Let $S=\{x_1,x_2,x_3,x_4\}\cup S'$. According to the construction of $S$, $S$ induces an $s$-clique in $G$. Since $|V(H_i)-\{x_1,x_2,x_3,x_4\}|$ is odd, $o(\overline{G}-S-P)\ge o(\overline{G}-P)-(s-4)$. By Observation \ref{ob1},
$$\chi(G-S)\ge \frac{(n-s)+o(\overline{G}-S-P)-|P|}2=\chi(G)-s+2=t+1.$$
This contradicts our assumption.\\

\noindent\textbf{Case 2.} $\omega(G)<s$.

As $\omega(G)\le s-1$ and $\alpha(G)=2$, the number of color classes of $G$ consisting of only one vertex is at most $s-1$. 
So \begin{equation}
n\ge(s-1)+2(\chi(G)-(s-1))=s+2t-1\ge3s-1.\label{eq2.2}
\end{equation} 
Since $\alpha(G)=2$, $\chi(G)\ge \frac n2. $
By Observation \ref{ob1}, 
\begin{equation}
|P|\le o(\overline{G}-P)\le\omega(G)\le s-1.\label{eq2.1}
\end{equation} 
Combining \eqref{eq2.2} and \eqref{eq2.1}, we have that $$n-|P|-o(\overline{G}-P)\ge2s-|P|.$$ 
Thus, we can pick a $2(s-\lceil|P|/2\rceil)$-element subset $S'\subseteq V(\overline{G}-P)$ that has an even number of vertices in common with each component of $\overline{G}-P$. Since each component of $\overline{G}-P$ is odd, $o(\overline{G}-P-S')\ge o(\overline{G}-P)$. Let $S=S'\cup P$. Then $$|S|=2s-2\lceil|P|/2\rceil+|P|\ge2s-1.$$
It follows that $\chi(G[S])\ge s$. By Observation \ref{ob1},
\begin{align*}
\chi(G-S)&\ge\frac{n-|S|+o(\overline{G}-S)}2\\
 &=\frac{n-|S|+o(\overline{G}-P)}2\\
&\ge\chi(G)-s+\lceil|P|/2\rceil\\
&= t+\lceil|P|/2\rceil-1.
\end{align*}
As $\chi(G-S)\le t$, we get that $|P|\le2$.

As each $H\in\mathcal{H}$ contains no $K_3$, we have that if $|H|\ge R(3,\ell)$, then $H$ contains an independent set of size at least $\ell$. We claim that there exists a component $H_0\subseteq\mathcal{H}$ such that $\alpha(H_0)\ge4$. Suppose to the contrary that for each $H\in\mathcal{H}$, $\alpha(H)\le3$.  Since $R(3,3)=6$ and $R(3,4)=9$, $o(\overline{G}-P)=k_1+k_3+k_5+k_7,$ 
\begin{equation}
k_1+3k_3+5k_5+7k_7=n-|P|\ge n-2\ge3s-3,\label{eq2.1.2}
\end{equation}
and
\begin{equation}
s-1\ge\omega(G)\ge k_1+2(k_3+k_5)+3k_7.\label{eq2.1.1}
\end{equation}
By \eqref{eq2.1.2} and \eqref{eq2.1.1}, we get that $k_1+3k_3+5k_5+7k_7\ge3(k_1+2(k_3+k_5)+3k_7)$. It implies that $2k_1+3k_3+k_5+2k_7\le0$, which is impossible as $n\ge3s-1\ge5$.

Let $x,y,z,w\in H_0$ form an independent set. Since $\overline{G}$ is triangle-free, $N_{\overline{G}}(v)$ is an independent set in $\overline{G}$ for each $v\in V(\overline{G})$. Let $F=N(x)\cup N(y)$. Then $|F|\le2\omega(G)\le2(s-1)$. Let $X$ be a set containing exactly one vertex from each component of $\mathcal{H}-H_0$ and $\{x,y,z\}\subseteq H_0$. By $\omega(G)\le s-1$, it is obvious that $|X|\le s-2$. Due to the fact that $n\ge3s-1$, $|V(\overline{G}-X)|\ge2s+1$. Next, we consider the following two cases.

\textbf{Subcase 2.1} $|P|\ge1$. That is, $|P|\in\{1,2\}$.

We get that $|V(\overline{G}-P-X)|\ge2s-1$. Note that if $|F|$ is odd, then $L_0=H_0-F-\{x,y,z\}\neq\emptyset$.  Now we start to construct $S'\subseteq V(\overline{G}-P-X)$ of size $2s-2$ by the following steps.

\begin{enumerate}
\item Put all vertices of $F$ into $S'$.
\item If $|F|$ is odd, then add a vertex $w\in L_0$ into $S'$. 
\item Add pairs from components of $\overline{G}-P-X$ until $|S'|=2s-2$.
\end{enumerate}

By the construction of $S'$, $S'$ has an even number of vertices in common with each component of $\overline{G}-P$. Since each component of $\overline{G}-P$ is odd and $x,y$ form two singleton components of $\overline{G}-S'-P$, $o(\overline{G}-P-S')\ge o(\overline{G}-P)+2$. Let $S=S'\cup P$. Then $$|S|=2s-2+|P|\ge2s-1.$$
It follows that $\chi(G[S])\ge s$. By Observation \ref{ob1},
\begin{align*}
\chi(G-S)&\ge\frac{n-|S|+o(\overline{G}-S)}2\\
 &=\frac{n-|P|+o(\overline{G}-P)}2+\frac{o(\overline{G}-S)-o(\overline{G}-P)+|P|-|S|}2\\
&\ge\chi(G)-s+\frac32> t.
\end{align*}
This contradicts our assumption.

\textbf{Subcase 2.2} $|P|=0$.

We claim that there are two vertices $u_1,u_2\in N_{\overline{G}}(x)$ such that for any $i\in[2]$, $u_i\notin N_{\overline{G}}(y)\cup N_{\overline{G}}(z)$. Suppose to the contrary that $N_{\overline{G}}(x)-u_1\subseteq N_{\overline{G}}(y)\cup N_{\overline{G}}(z)$. Let $F'=N_{\overline{G}}(x)\cup N_{\overline{G}}(y)\cup N_{\overline{G}}(z)$. It follows that $|F'|\le2s-1.$
Since $w\notin F'$, $L_0=H_0-F'-\{x,y,z\}\neq\emptyset$. If $|F'|$ is even, then $|L_0|\ge2$. Now we start to construct $S\subseteq V(\overline{G}-X)$ of size $2s-1$ by the following steps.
\begin{enumerate}
\item Put all vertices of $F'$ into $S$.
\item If $|F'|$ is even, then add $w\in L_0$ into $S$. 
\item Add pairs from components of $\overline{G}-X$ until $|S|=2s-1$.
\end{enumerate}
By the construction of $S$, $S$ has an even number of vertices in common with each component of $\overline{G}$ except $H_0$. Since each component of $\overline{G}$ is odd and $x,y,z$ form three singleton components of $\overline{G}-S$, $o(\overline{G}-S)\ge o(\overline{G})+3$. Note that $\chi(G[S])\ge s$. By Observation \ref{ob1},
\begin{align*}
\chi(G-S)&\ge\frac{n-|S|+o(\overline{G}-S)}2\\
 &=\frac{n+o(\overline{G})}2+\frac{o(\overline{G}-S)-o(\overline{G})-|S|}2\\
&\ge\chi(G)-s+2\ge t+1.
\end{align*}
This contradicts our assumption.

Let $U_1=\{x\}\cup N_{\overline{G}}(x)\cup N_{\overline{G}}(u_1)\cup N_{\overline{G}}(u_2)$ and $U_2$ be a set containing exactly one vertex from each component of $\mathcal{H}-H_0$. Clearly, $U_1\subset H_0$ and $U_2$ is independent. Let $Y=U_1\cup U_2$. Next, we will choose a set $S'\subseteq V(\overline{G})-Y$ of size $2s-3-|N_{\overline{G}}(x)|$. Firstly, let us calculate the size of $V(\overline{G})-Y$. Note that $$\chi(G)=\frac{n+o(\overline{G})}2.$$ We have that $$n=2s+2t-2-o(\overline{G}).$$
As $N_{\overline{G}}(v)$ is an independent set for any $v\in V(G)$, we get that $U_2\cup N_{\overline{G}}(u)$ is independent for any $u\in H_0$. So for any $u\in H_0$,

\begin{equation*}\label{eq2.2.1}
 |U_2\cup N_{\overline{G}}(u)|=o(\overline{G})+|N_{\overline{G}}(u)|-1\le \omega(G)\le s-1.
\end{equation*}

It follows that
\begin{align*}\label{eq2.2.2}
 |V(\overline{G})-Y|&=2s+2t-2-o(\overline{G})-|U_2\cup N_{\overline{G}}(u_1)\cup N_{\overline{G}}(u_2)|-|N_{\overline{G}}(x)| \\
 &\ge 2s+2t-2-(o(\overline{G})+|N_{\overline{G}}(u_1)|)-(|U_2\cup N_{\overline{G}}(u)|)-|N_{\overline{G}}(x)|\\
 &\ge 2s+2t-2-s-(s-1)-|N_{\overline{G}}(x)|\\
 &=2t-1-|N_{\overline{G}}(x)|>2s-3-|N_{\overline{G}}(x)|.
\end{align*}

Since $y,z\notin N_{\overline{G}}(x)$, $L_0=H_0-Y\neq\emptyset$. If $|N_{\overline{G}}(x)|$ is even, then $|L_0|\ge2$. Now we start to construct $S'\subseteq V(\overline{G}-Y)$ by the following steps.
\begin{enumerate}
 \item If $|N_{\overline{G}}(x)|$ is even, then add $y\in L_0$ into $S'$. 
\item Add pairs from components of $\overline{G}-Y$ until $|S'|=2s-3-|N_{\overline{G}}(x)|$.
\end{enumerate}
Let $S=S'\cup N_{\overline{G}}(x)$. Then $|S|=2s-3$. As $\alpha(G)=2$, so $\chi(S-u_1-u_2)\ge s-2.$ In $\overline{G}$, note that $u_i$ has no neighbors in $S$ for $i\in[2]$. That is, in graph $G$, $u_i$ is adjacent to all vertices in $S$. It follows that $\chi(S)\ge s-2+2=s$.

By the construction of $S$, $S$ has an even number of vertices in common with each component of $\overline{G}$ except $H_0$. Since each component of $\overline{G}$ is odd and $x$ becomes a singleton component of $\overline{G}-S$, $o(\overline{G}-S)\ge o(\overline{G})+1$. By Observation \ref{ob1},
\begin{align*}
\chi(G-S)&\ge\frac{n-|S|+o(\overline{G}-S)}2\\
 &=\frac{n+o(\overline{G})}2+\frac{o(\overline{G}-S)-o(\overline{G})-|S|}2\\
&\ge\chi(G)-s+2\ge t+1.
\end{align*}
This contradicts our assumption. Thus, a counterexample $G$ does not exist for the above statement.
\end{proof}


\bigskip






\begin{thebibliography}{99}

\bibitem{BKPS09} J. Balogh, V. Kostochka, N. Prince, and M. Stiebitz, The Erd\H{o}s-Lo\'{v}asz Tihany conjecture for quasi-line graphs. Discrete Math. 309 (2009), 3985-3991.

\bibitem{BJ99} W. G. Brown and H. A. Jung, On odd circuits in chromatic graphs, Acta Math. Acad. Sci. Hungar. 20 (1999), 129-134.

\bibitem{EL68} P. Erd\H{o}s, Problem 2, In: Theory of Graphs (P. Erd\H{o}s and G. Katona, Eds.), Proc. Colloq. Tihany, Hungary, September 1966, Academic Press, New York, 1968, p. 361.

\bibitem{KS08} A. V. Kostochka, M. Stiebitz, Partitions and edge colorings of multigraphs, Electron. J. Combin., 15 (2008), N25.

\bibitem {M86} N. N. Mozhan, On doubly critical graphs with chromatic number five, Technical Report 14, Omsk Institute of Technology, 1986 (in Russian).

\bibitem{S87a} M. Stiebitz, $K_5$ is the only double-critical 5-chromatic graph, Discrete Math. 64 (1987), 91-93.

\bibitem{S87b} M. Stiebitz, On k-critical n-chromatic graphs. In: Colloquia Mathematica Soc. J\'{a}nos Bolyai
52, Combinatorics, Eger (Hungary), 1987, 509-514.

\bibitem{WY} Y. Wang, G. Yu, Enhancing the Erd\H{o}s-Lov\'asz Tihany Conjecture for line graphs of multigraphs (Submitted).

\end{thebibliography}
\end{document}